\documentclass[12pt,a4paper,twoside,leqno]{article}
\usepackage[english]{babel}
\usepackage{amssymb}
\usepackage{inputenc}
\newcommand{\ba}{\begin{array}}
\newcommand{\ea}{\end{array}}
\usepackage{amssymb}
\newtheorem{theorem}{Theorem}[section]

\newtheorem{definition}[theorem]{Definition}
\newtheorem{proposition}[theorem]{Proposition}
\newtheorem{example}[theorem]{Example}
\newtheorem{lemma}[theorem]{Lemma}
\newtheorem{corollary}[theorem]{Corollary}

\usepackage{cite}
\DeclareSymbolFont{rsfs}{U}{rsfs}{m}{n}
\DeclareSymbolFontAlphabet{\rsfsmathscr}{rsfs}
\usepackage[mathscr]{eucal}
\usepackage{amscd}
\usepackage{floatflt}
\usepackage{amsmath}
\usepackage{amssymb}

\def\e{{\bf 1}\!\!{\rm I}}
\numberwithin{equation}{section}
\begin{document}

\author{S. Albeverio $^{1},$ Sh. A. Ayupov $^{2,  *},$ A. A. Rakhimov, $^3$ R. A. Dadakhodjaev $^{4}$}

\title{\bf  Index for Finite Real Factors}

\date{}

\maketitle

\begin{abstract}
For real factors the notions of the coupling constant and the index are introduced and
investigated. The possible values of the index for type II$_1$ real factors are
calculated, in a similar way as for the complex case.

\end{abstract}
\medskip
\medskip

$^1$ Institut f\"{u}r Angewandte Mathematik, Universit\"{a}t Bonn,
Wegelerstr. 6, D-53115 Bonn (Germany); SFB 611; HCM;  BiBoS; IZKS; CERFIM
(Locarno); Acc. Arch. (USI),

e-mail address:\emph{albeverio@uni-bonn.de}

\medskip
$^2$ Institute of Mathematics and information  technologies,
Uzbekistan Academy of Sciences, Dormon Yoli str. 29, 100125,
Tashkent (Uzbekistan),

e-mail address: \emph{sh\_ayupov@mail.ru}
\medskip

$^{3}$ The Tashkent institute of railways and engineers,
Adilkhodjaev str. 1, Tashkent(Uzbekistan),

e-mail address: \emph{rakhimov@ktu.edu.tr}
\medskip

$^4$ Institute of Mathematics and information  technologies,
Uzbekistan Academy of Sciences, Dormon Yoli str. 29, 100125,
Tashkent (Uzbekistan),

e-mail address: \emph{rashidkhon@mail.ru}

\medskip



\medskip
\medskip

* Corresponding author

\section{{\large Introduction} }

In 1930's Von Neumann and Murray introduced the notion of coupling
constant for finite factors (see \cite{MN1,MN2,MN3}). In 1983, V.Jones
suggested a new approach to this notion, defined the notion of index for type II$_1$ factors,
and proved a surprising theorem on values of the index for
subfactors (see \cite{J}). He also introduced a very important
technique in the proof of this theorem: the towers of algebras.
Since then this theory has become a focus of many fields in mathematics and
physics (\cite{J1}). In \cite{K1}, H.Kosaki extended the notion of
the index to an arbitrary (normal faithful) expectation from a factor onto a
subfactor. While Jones' definition of the index is based on the
coupling constant, Kosaki's definition of the index of an
expectation relies on the notion of spatial derivatives due to
A.Connes \cite{C} as well as the theory of operator-value weights
due to U.Haagerup \cite{H}. In \cite{K1,K2}, it was shown that many
fundamental properties of the Jones index in the type II$_1$ case
can be extended to the general setting. At the present time, the theory
of index thanks to works by V.Jones, P.Loi, R.Longo, H.Kosaki and
other mathematicians is deeply developed and has many applications
in the theory of operator algebras and physics (see also \cite{PL,RL}).\\[-2mm]

Unlike to the complex case for real factors the notion of coupling
constant (therefore the notion of index as well) has not been
investigated. In the present paper the notions of the real coupling
constant and the index for finite real factors are introduced and
investigated. The main tool in our approach is the reduction of real factors to
involutive *-anti-automorphisms of their complex enveloping von Neumann algebras.

\section{{\large Preliminaries} }

Let $B(H)$ be the algebra of all bounded linear operators on a
complex Hilbert space $H$. A weakly closed *-subalgebra
$\mathfrak{A}$ containing the identity operator $\e$ in $B(H)$ is called a
W$^*$-{\it algebra}. A real *-subalgebra $\Re\subset B(H)$ is called
a {\it real} W$^*$-{\it algebra} if it is closed in the weak
operator topology and \ $\Re \cap i\Re = \{0\}$. A real
W$^*$-algebra $\Re$ is called a {\it real factor} if its center
$Z(\Re)$ consists of the elements $\{\lambda\e , \lambda\in\mathbb{R}\}$.
We say that a real W$^*$-algebra $\Re$ is
of the type I$_{fin}$, I$_{\infty}$, II$_1,$ II$_{\infty},$ \ or \
III$_{\lambda}$, $(0 \leq \lambda \leq 1)$ \ if the enveloping
W$^*$-algebra $\mathfrak{A}(\Re)$ has the corresponding type in the
ordinary classification of W$^*$-algebras. A linear mapping $\alpha$
of an algebra into itself with $\alpha(x^*)=\alpha(x)^*$ \ is called
an *-{\it automorphism} if $\alpha(xy)=\alpha(x)\alpha(y)$;
it is called an {\it involutive *-antiautomorphism} if $\alpha(xy)=\alpha(y)\alpha(x)$
and $\alpha ^2(x)=x$. If $\alpha$ is an involutive
*-antiautomorphism of a W$^*$-algebra $M$, we denote by $(M,\alpha)$
the real W$^*$-algebra generated by $\alpha$, i.e. $(M,\alpha)=\{x\in M
: \ \alpha(x)=x^*\}$. Conversely, every real W*-algebra $\Re$ is of
the form $(M,\alpha)$, where $M$ is the complex envelope of $\Re$ and
$\alpha$ is an involutive *-antiautomorphism of $M$ (see \cite{ARU,G,R}).
Therefore we shall identify from now on the real von Neumann algebra $\Re$
with the pair $(M,\alpha)$.


\section{{\large Canonical representation} }

Let $M$ ($\subset B(H)$) be a finite factor and let $\tau$ be the
unique faithful normal tracial state of $M$. If $\alpha$ is an
involutive *-antiautomorphism of $M$, then it is clear that $\tau$
is automatically $\alpha$-invariant. Denote by $L^2(M)$ the completion of
$M$ with respect to the norm $\|x\|_2=\tau(x^*x)^{1/2}$.
Similarly by $L^2(M,\alpha)$ we denote the completion of the real
factor $(M,\alpha)$. Then it is obvious that the Hilbert space $L^2(M)$
and the algebra $B(L^2(M))$ of all bounded linear operators on it
are the complexifications of the real Hilbert space $L^2(M,\alpha)$ and
of $B_r(L^2(M,\alpha))$, respectively, where $B_r(L^2(M,\alpha))$ is
the algebra of all bounded linear operators on the real Hilbert
space $L^2(M,\alpha)$. Moreover, it is easy to show that the Hilbert spaces
$L^2(M,\alpha)$ and $L^2(M)$ are separable.

\vspace{0.5cm}

For each $x\in M$, set $\lambda(x)y=xy$, for all $y\in M$. Clearly,
$\|\lambda(x)y\|_2\leq\|x\|\|y\|_2$. Thus $\lambda$ can be uniquely
extended to a bounded linear ope\-rator on $L^2(M)$, still denoted
by $\lambda(x)$. Then we obtain a faithful W$^*$-representation
$(\lambda , L^2(M))$ of $M$. In a similar way, taking the map
$\lambda _r$ defined as $\lambda _r(x)y=xy$ (for all $x,y\in
(M,\alpha)$) we obtain a faithful real $*$-representation
$(\lambda _r, L^2(M,\alpha))$ of $(M,\alpha)$.

\begin{theorem} \label{teorem002} The map $\beta :\lambda(M)\to\lambda(M)$
defined as \ $\beta(\lambda _x) = \lambda _{\alpha(x)}$
is an involutive *-antiautomorphism of $\lambda(M)$. Moreover,
$\beta$ and $\alpha$ are also related in the following way:
\ $(M,\alpha)_\beta = \lambda _r(M,\alpha)$, where
$(M,\alpha)_\beta=\{\lambda _x\in \lambda(M) : \beta(\lambda
_x)=\lambda _x^*\}$ \ is the real W$^*$-algebra, generated by $\beta$,
i.e. $(M,\alpha)_\beta=(\lambda(M),\beta)$.
\end{theorem}

\noindent {\bf Proof}. The first part of the assertion is trivial.
Further, let $\lambda _x\in (M,\alpha)_\beta$. Since $\beta(\lambda _x)=\lambda _x^*$, then
$\lambda _{\alpha(x)}=\lambda _{x^*}$. Hence $\alpha(x)=x^*$, i.e.
$x\in (M,\alpha)$. Then from
$$
\lambda _x\in \lambda(M)\subset B(L^2(M)) = B_r(L^2(M,\alpha)) + i B_r(L^2(M,\alpha))
$$
we have \ $(M,\alpha)_\beta \subset B_r(L^2(M,\alpha))$. Hence \
$(M,\alpha)_\beta \subset \lambda _r(M,\alpha)$,
since $\lambda _r(M,\alpha)=\{\lambda ^r_x\in B_r(L^2(M,\alpha))$ : for
$\alpha(x)=x^*$ and $\alpha(y)=y^*$, $\lambda ^r_x(y):=xy\}$.\\[-2mm]

Now let \ $\lambda ^r_x\in\lambda _r(M,\alpha)$. Then
$\alpha(x)=x^*$ and $\lambda ^r_x\in\lambda _r(M,\alpha)\subset
\lambda(M)$. Hence $\beta(\lambda ^r_x) = \lambda ^r_{\alpha(x)} =
\lambda ^r_{x^*} = (\lambda ^r_x)^*$, therefore \ $\lambda ^r_x \in
(M,\alpha)_\beta$. \ $\Box$

\begin{corollary} \label{corol001} $\lambda _r(M,\alpha)$ is a real W$^*$-algebra, and
$\lambda(M)$ is the complexification of $\lambda _r(M,\alpha)$, i.e.
$\lambda _r(M,\alpha) + i\lambda _r(M,\alpha) = \lambda(M)$.
Moreover, \ $\{\lambda _r, L^2(M,\alpha)\}$ is a faithful real
W$^*$-representation of $(M,\alpha)$.
\end{corollary}

This representation will be called the canonical W*-representation of $(M,\alpha)$.

\section{{\large Commutant of the canonical representation} }


Since $\|x\|_2 = \|x^*\|_2$ for all $x\in M$, the map $J:x\to x^*$
can be uniquely extended to a conjugate linear isometry on $L^2(M)$,
still denoted by $J$. From the theory of W$^*$-algebras it is
well-known that $\lambda(M)' = J\lambda(M)J$ and $\lambda(M) =
J\lambda(M)'J$. Similarly to Theorem \ref{teorem002} and Corollary
\ref{corol001} we have the following assertion

\begin{theorem} \label{teorem0030}
The map $\beta ' :\lambda(M)'\to \lambda(M)'$ defined as $\beta '(
\cdot ) = J \beta(J \cdot J) J$, is an involutive *-antiautomorphism
of $\lambda(M)'$. The set \ $\lambda _r(M,\alpha)'=\{\lambda _{x'}\in
\lambda(M)' : \quad \beta '(\lambda _{x'}) = \lambda _{x'}^*\}$ \
is a real W$^*$-algebra, and $\lambda(M)'$ is the complexification of $\lambda
_r(M,\alpha)'$, i.e. $\lambda _r(M,\alpha)'+i\lambda _r(M,\alpha)'=\lambda(M)'$.
\end{theorem}

We have the following connection between
$\lambda _r(M,\alpha)$ and $\lambda _r(M,\alpha)'$.

\begin{theorem} \label{teorem003}
\quad $\lambda _r(M,\alpha)' = J \lambda _r(M,\alpha) J$ .
\end{theorem}

\noindent {\bf Proof}. Since $\lambda _x\in\lambda _r(M,\alpha)$ implies that
$J\lambda _xJ\in J\lambda _r(M,\alpha)J$ \ and \ $\beta(\lambda _x) = \lambda _x^*$,
we have
$$
\beta '(J\lambda
_xJ)=J\beta(JJ\lambda _xJJ)J=J\beta(\lambda _x)J=J\lambda
_x^*J=(J\lambda _xJ)^* .
$$
Hence \ $J\lambda _xJ\in \lambda _r(M,\alpha)'$, \ i.e. \ $J\lambda _r(M,\alpha)J \subset \lambda _r(M,\alpha)'$.\\[-3mm]

Conversely, let $\lambda _{x'}\in \lambda _r(M,\alpha)' \subset \lambda(M)'$ =$J\lambda(M)J$. Then
$\lambda _{x'}=J\lambda _yJ$, for some $\lambda _y\in \lambda(M)$.
Since \ $\beta(\lambda _{x'})=\lambda _{x'}^*$, we have \ $\beta '(J\lambda _yJ)=J\lambda _y^*J$, \ i.e. \
$J\beta(JJ\lambda _yJJ)J=J\lambda _y^*J$.
Hence
\ $J^2\beta(\lambda _y)J^2=J^2\lambda _y^*J^2$, \ i.e. \ $\beta(\lambda _y)=\lambda _y^*$.
\ Therefore \ $\lambda _y\in \lambda _r(M,\alpha)$. Thus we obtain
$\lambda _{x'}=J\lambda _yJ=J\lambda _r(M,\alpha)J$, and
therefore \ $\lambda _r(M,\alpha)' \subset J\lambda _r(M,\alpha)J$. \ $\Box$

\begin{theorem} \label{teorem0031} The real W$^*$-algebra $\lambda _r(M,\alpha)'$ is
the commutant of $\lambda _r(M,\alpha)$ in the algebra $B_r(L^2(M,\alpha))$, i.e.
$\lambda _r(M,\alpha)'=\{\lambda _x\in B_r(L^2(M,\alpha)) : \lambda
_x\lambda _y=\lambda _y\lambda _x , \ \forall \ \lambda _y\in\lambda
_r(M,\alpha)\}$
\end{theorem}

\noindent {\bf Proof}. Similarly to the proof of Theorem
\ref{teorem002} for $\beta '(\lambda _x)=\lambda _x^*$ we have
$\lambda _x\in B_r(L^2(M,\alpha))$. Therefore $\lambda
_r(M,\alpha)'\subset B_r(L^2(M,\alpha))$. On the other hand for any
$\lambda _x\in \lambda _r(M,\alpha)'$ $\subset \lambda(M)'$ and
$\lambda _y\in \lambda _r(M,\alpha)$ $\subset\lambda(M)$,
\ we have \ $\lambda _x\lambda _y=\lambda _y\lambda _x$. \ $\Box$

\section{{\large Relations between faithful nondegenerate W$^*$-represen\-ta\-tions and the canonical representation} }


\begin{theorem} \label{teorem004}
Let $M_1\subset B(H_1)$ and $M_2\subset B(H_2)$ be two
$W^*$-algebras and let $\alpha _i$ be an involutive *-antiautomorphism
of $M_i$, $i=1,2$. If \ $\Phi : M_1 \to M_2$ is a normal
*-homomorphism with $\Phi \circ \alpha _1 = \alpha _2 \circ \Phi$,
then
$$
\Phi = \Phi _3 \circ \Phi _2 \circ \Phi _1 ,
$$
where \\[-3mm]

$\Phi _1$ is a *-homomorphism from $M_1$ \ onto \
$M_1{\overline{\otimes}} \ \mathbb{C}\e _L$ \ with \
$\Phi _1 \circ \alpha _1 = {\tilde{\alpha _1}} \circ \Phi _1$ \ defined as
\ $\Phi(a)=a\otimes \e _L$, \ where $\e _L$ is the identity operator on an appropriate
Hilbert space $L$ \ and \ ${\tilde{\alpha _1}}=\alpha _1 \otimes id$; \\[-3mm]

$\Phi _2$ is a *-homomorphism from $M_1{\overline{\otimes}} \ \mathbb{C}\e
_L$ \ onto \ $(M_1 {\overline{\otimes}} \ \mathbb{C}\e _L)p'$ \ with \
$\Phi _2 \circ {\tilde{\alpha _1}} = {\overline{\alpha}}_1 \circ \Phi _2$
\ defined as \ $\Phi _2(a\otimes \e _L)=(a\otimes \e _L)p'$, \ where
$p'$ is a projection from $(M_1{\overline{\otimes}} \ \mathbb{C}\e _L)'$ with
${\tilde{\alpha _1}}'(p') = p'$ \ and \
${\tilde{\alpha _1}}'=J_1{\tilde{\alpha _1}}(J_1( . )J_1)J_1
\otimes id$, \ ${\overline{\alpha}}_1( \cdot \ p') = {\tilde{\alpha _1}}( \cdot ) p'$; \\[-3mm]

$\Phi _3$ is a *-isomorphism from $(M_1 {\overline{\otimes}} \ \mathbb{C}\e _L)p'$ \ to \
$M_2$ \ with \ $\Phi _3 \circ {\tilde{\alpha _1}} = \alpha _2 \circ \Phi _3$.
\end{theorem}

\noindent {\bf Proof}. First we assume that $(M_2, \alpha _2)$
admits a cyclic vector $\eta$. In this case \ ${\overline{(M_2,
\alpha _2) \eta}} = H_2^r$ \ is a real Hilbert space and
$$
{\overline{M_2 \eta}} = {\overline{(M_2, \alpha _2) \eta}} + i
{\overline{(M_2, \alpha _2) \eta}} = H_2^r + i H_2^r = H_2 ,
$$
hence $\eta$ is a cyclic vector of $M_2$. Since \ $\Phi \circ \alpha
_1 = \alpha _2 \circ \Phi$, \ for all $a\in (M_1, \alpha _1)$
\ we have \ $\alpha _2\bigl( \Phi(a)\bigr) = \Phi\bigl( \alpha
_1(a)\bigr) = \Phi(a^*) = \Phi(a)^*$, i.e. \ $\Phi(a)\in (M_2,
\alpha _2)$. Hence \ $\Phi\bigl( (M_1, \alpha _1) \bigr) \subset
(M_2, \alpha _2)$. \ Define a functional $\varphi$ by
$$
\varphi(a) = <\Phi(a)\eta , \eta > , \quad \ a\in (M_1, \alpha _1) .
$$
Obviously, $\varphi$ is a normal positive functional on $(M_1,
\alpha _1)$. We can extend $\varphi$ by linearity to a functional on $M_1$
(still denoted by $\varphi$) such that
$$
\varphi(a+ib) = \varphi(a)+i\varphi(b) , \qquad a,b\in (M_1, \alpha _1) ,
$$
which clearly also is a normal positive functional.
Let $H_1^r$ be a real Hilbert space with $H_1^r+iH_1^r=H_1$ such that
$(M_1, \alpha _1)\subset B(H_1^r)$.
By \cite[4.2.1]{Li2} there is a sequence $(\xi _n)\subset H_1^r$ with
$\sum\limits _n\|\xi _n\|^2<\infty$ such that
$\varphi(a)=\sum\limits _n<a\xi _n,\xi _n>$, for all $a\in (M_1,
\alpha _1)$. Set $L_r=\ell _2^r=\{(x_n)\subset \mathbb{R} : \sum
_nx_n^2<\infty\}$, $L=L_r+iL_r$, $\xi=(\xi _n)\subset H_1^r\otimes
L_r$ and $\Phi _1(a)=a\otimes \e _L$ for all $a\in M_1$. Then \
$\Phi _1$ is a map from $M_1$ to $M_1{\overline{\otimes}} \ \mathbb{C}\e _L$ \ and
\begin{eqnarray*}
(\Phi _1 \circ \alpha _1)(a) & = & \Phi _1(\alpha _1(a)) = \alpha _1(a)\otimes \e _L = (\alpha _1 \otimes id)(a\otimes \e _L) \\[1mm]
                             & = & {\tilde{\alpha _1}}(\Phi _1(a)) = ({\tilde{\alpha _1}} \circ \Phi _1)(a),
\end{eqnarray*}
i.e. \ $\Phi _1 \circ \alpha _1 = {\tilde{\alpha _1}} \circ \Phi
_1$. Moreover, for all $a\in (M_1, \alpha _1)$ \ we have
$$
<\Phi _1(a)\xi , \xi > = <(a\otimes \e _{L_r})\xi , \xi > = \sum _n
<a\xi _n, \xi _n> = \varphi(a) .
$$
Let $p'$ be the projection from $H_1^r\otimes L_r$ to
${\overline{\Phi _1\bigl( (M_1, \alpha _1) \bigr)\xi}}$. Then for all
$x=a\otimes \e _{L_r}\in ((M_1, \alpha _1) {\overline{\otimes}} \
\mathbb{R}\e _{L_r})$ \ we have
\begin{eqnarray*}
(p'x)\xi & = & p'((a\otimes \e _{L_r})\xi) = p'(\Phi _1(a)\xi) = \Phi _1(a)\xi \\[1mm]
         & = & (a\otimes \e _{L_r})\xi = x\xi = x((\e\otimes \e _{L_r})\xi) = x(\Phi _1(\e)\xi) \\[1mm]
         & = & x(p'(\Phi _1(\e)\xi)) = x(p'(\xi)) = (xp')\xi .
\end{eqnarray*}
Similarly, for all $\gamma\in H_1^r\otimes L_r$ with \ $\gamma
\not = \xi$ \ we also obtain
\begin{eqnarray*}
(p'x)\gamma & = & p'(\Phi _1(a)\gamma) = \theta = x(\theta) = x(p'(\Phi _1(\e)\gamma))\\[1mm]
            & = & xp'((\e\otimes \e _{L_r})\gamma) = xp'(\gamma) .
\end{eqnarray*}
Therefore $p'x = xp'$, i.e. \ $p'\in ((M_1, \alpha _1)
{\overline{\otimes}} \ \mathbb{R}\e _{L_r})'$. Hence \ $p'\in (M_1
{\overline{\otimes}} \ \mathbb{C}\e _L)'$ and for ${\tilde{\alpha
_1}}' = J_1{\tilde{\alpha _1}}(J_1( . )J_1)J_1 \otimes id$ \ we have
\ ${\tilde{\alpha _1}}'(p') = p'$.\\[-2mm]

\noindent Define the map \ $\Phi _2 : M_1{\overline{\otimes}} \
\mathbb{C}\e _L \to (M_1 {\overline{\otimes}} \ \mathbb{C}\e _L)p'$
\ as \ $\Phi _2(a\otimes \e _L)=(a\otimes \e _L)p'$, \ $a\in M_1$.
Then
\begin{eqnarray*}
(\Phi _2 \circ {\tilde{\alpha _1}})(a\otimes \e _L) & = & \Phi
_2({\tilde{\alpha _1}}(a\otimes \e _L)) =
\Phi _2(\alpha _1(a)\otimes \e _L)  \\[2mm]
& = & (\alpha _1(a)\otimes \e _L)p' = {\tilde{\alpha _1}}(a\otimes \e _L)p'\\[2mm]
& = & {\overline{\alpha}}_1\bigl( (a\otimes \e _L)p'\bigr) = {\overline{\alpha}}_1(\Phi _2(a\otimes \e _L))\\[2mm]
& = & ({\overline{\alpha}}_1 \circ \Phi _2)(a\otimes \e _L) ,
\end{eqnarray*}
hence \ $\Phi _2 \circ {\tilde{\alpha _1}} = {\overline{\alpha}}_1
\circ \Phi _2$. Since \ $p'\xi = p'((\e\otimes \e _L)\xi) = p'(\Phi
_1(\e)\xi) = \Phi _1(\e)\xi = \xi$, we have
\begin{eqnarray*}
<(\Phi _2\circ \Phi _1)(a)\xi , \xi > & = & <(\Phi _2(a\otimes \e _{L_r}))\xi , \xi > = <(a\otimes \e _{L_r})p'\xi , \xi >\\[2mm]
& = & <(a\otimes \e _{L_r})\xi , \xi > = <\Phi _1(a)\xi , \xi > =
\varphi(a) ,
\end{eqnarray*}
for all \ $a\in (M_1, \alpha _1)$, i.e. \ $\varphi(a)=<(\Phi _2\circ
\Phi _1)(a)\xi ,\xi >$.

\vspace{0.3cm}

\noindent Now, define a linear map \ $u : \Phi\bigl( (M_1, \alpha
_1) \bigr)\eta \to p'(H_1^r\otimes L_r)$ \ as follows:
$$
u \Phi(a)\eta = (\Phi _2\circ \Phi _1)(a)\xi = p'(a\xi _n) = (a\xi
_n) \quad (a\in (M_1, \alpha _1)).
$$
Since \ $u \Phi(a)\eta = (\Phi _2\circ \Phi _1)(a)\xi$ \ and \
$<\Phi(a)\eta ,\eta > = \varphi(a) = <(\Phi _2\circ \Phi _1)(a)\xi
,\xi >$ ($a\in (M_1, \alpha _1)$), it follows that \
$\|u\Phi(a)\eta\|' = \|\Phi(a)\eta\|_2^r$, i.e. the map $u$ is an isometry,
where \ $\|\cdot \|_2^r$ is the norm of the space $H_2$ and $\|\cdot \|'$
is the norm of the space \ $H_1^r\otimes L_r$. Moreover, since \
$\Phi\bigl( (M_1, \alpha _1) \bigr)\eta = (M_2, \alpha _2)\eta$, \
$(\Phi _2\circ \Phi _1)\bigl( (M_1, \alpha _1) \bigr)\xi = \Phi _1\bigl( (M_1, \alpha _1) \bigr)\xi$, \ and\\[-1mm]

\qquad ${\overline{\Phi\bigl( (M_1, \alpha _1) \bigr)\eta}} = {\overline{(M_2, \alpha _2) \eta}} = H_2^r$ ,\\

\qquad ${\overline{(\Phi _2\circ \Phi _1)\bigl( (M_1, \alpha _1)
\bigr)\xi}} =
{\overline{\Phi _1\bigl( (M_1, \alpha _1) \bigr)\xi}} = p'(H_1^r\otimes L_r)$,\\

\noindent $u$ can be extended to a unitary operator ${\overline{u}}
: H_2^r \to p'(H_1^r\otimes L)$. Clearly,
\begin{eqnarray} \label{eq-for001}
{\overline{u}} \Phi(a) {\overline{u}}^{\ -1} = \Phi _2\circ \Phi
_1(a) , \quad a\in (M_1, \alpha _1) .
\end{eqnarray}
Therefore we can define a spatial real *-isomorphism $\Phi _3 : \bigl(
(M_1, \alpha _1) {\overline{\otimes}} \ \mathbb{R}\e _{L_r}\bigr)p'$
$\to (M_2, \alpha _2)$ \ as \ $\Phi _3( . ) = {\overline{u}}^{\ -1}
( . ) {\overline{u}}$, and it can be extended to a spatial
*-isomorphism (still denoted by $\Phi _3$) \ $\Phi _3 : (M_1
{\overline{\otimes}} \ \mathbb{C}\e _L)p' \to M_2$
\ as \ $\Phi _3(a+ib)=\Phi _3(a)+i\Phi _3(b)$, \ where \ $a,b\in \bigl( (M_1, \alpha _1) {\overline{\otimes}} \ \mathbb{R}\e _{L_r}\bigr)p'$.
Then, by \eqref{eq-for001} we have \ $\Phi = \Phi _3 \circ \Phi _2 \circ \Phi _1$.

\vspace{0.3cm}

Considering now the general case, the real Hilbert space $H_2^r$ with \ $H_2^r+iH_2^r=H_2$ \ can be decomposed as
$$
H_2^r = \oplus _l H^l_2 \quad {\rm and} \quad H^l_2 =
{\overline{(M_2, \alpha _2) \eta _l}}, \quad {\rm where} \quad \eta _l\in H_2^r ,
\quad {\rm for} \ l\in \mathbb{N} .
$$
Let $q_l' : H_2^r\to {\overline{(M_2, \alpha _2) \eta _l}}=H^l_2$ \
be the projection. Then \ $q_l'\in (M_2, \alpha _2)'$, \ for all \
$l$. For each $l$, \ $\Phi _l=q_l'\Phi : (M_1, \alpha _1)\to (M_2,
\alpha _2)q_l'$ \ is a normal *-homomorphism, which can be extended to
a normal *-homomorphism \ $\Phi _l : M_1\to M_2q_l'$. Then, by the
above argument \ $\Phi _l = \Phi _3^{(l)} \circ \Phi _2^{(l)} \circ
\Phi _1^{(l)}$, \ for all \ $l$. Set \ $\Phi _i = \oplus _l \Phi
_i^{(l)}$, \ $i=1,2,3$. Then \ $\Phi = \Phi _3 \circ \Phi _2 \circ
\Phi _1$ \ and the maps \ $\Phi _3, \Phi _2, \Phi _1$ satisfy all our
conditions. \ $\Box$\\[2mm]

\begin{theorem} \label{teorem005}
Let $M$ be a finite factor and let $\alpha$ be an involutive
*-anti\-automorphism of $M$. If \ $\{\pi , H\}$ \ is a faithful
nondegenerate W*-repre\-sentation of $M$ and \ $\pi \circ \alpha =
{\tilde{\alpha}} \circ \pi$ \ for an involutive *-antiautomorphism
${\tilde{\alpha}}$ of $\pi(M)$, then there exist a projection \ $p'\in
(\lambda _r(M,\alpha)\otimes \e _{K_r})'$, and a unitary operator $u
: H_r\to p'(L^2(M,\alpha)$ $\otimes K_r)$ such that
$$
u\pi(x) = (\lambda(x)\otimes \e _K)u , \qquad \ x\in M ,
$$
i.e., the real W*-algebras $\pi(M,\alpha)$
$(=(\pi(M),{\tilde{\alpha}}))$ and $(\lambda _r(M,\alpha) \otimes \e
_{K_r})p'$ \ are spatially *-isomorphic and therefore the W*-algebras
$\pi(M)$ and $(\lambda(M) \otimes \e _K)p'$ \ are also spatially
*-isomorphic; where $K_r$ is a separable infinite dimensional
Hilbert space, and $K=K_r+iK_r$.
\end{theorem}

\noindent {\bf Proof}. \ Set $M_1=\lambda(M)$ and $M_2=\pi(M)$.
Define the map \ $\Phi : M_1\to M_2$ \ by \ $\Phi(\lambda(x))=\pi(x)$.
Then $\Phi$ is a *-isomorphism and \ $\Phi\bigl(\lambda
_r(M,\alpha)\bigr) \subset (\pi(M),{\tilde{\alpha}})$.
Now the conclusion follows immediately from Theorem \ref{teorem004} and the separability of $H$. \ $\Box$\\


\section{{\large The coupling constants for real factors} }

If $M$ ($\subset B(H)$) is a finite factor with the finite commutant
$M'$, the {\it coupling constant} $\dim _M(H)$ of $M$ is defined as
tr$_M(E^{M'}_\xi)/$tr$_{M'}(E^M_\xi)$, where $\xi$ is a non-zero
vector in $H$, tr$_A$ denotes the normalized trace and $E^A_\xi$ is
the projection onto the closure of the subspace $A\xi$. This
definition, due to Murray and von Neumann in \cite{MN1}, is
independent of $\xi$. We recall some properties of the coupling
constant (\cite{MN1,MN2}, see also \cite{J},\cite{D},\cite[Ch. 17]{Li1})
{\small
\begin{subequations} \label{eq-for01}
\begin{eqnarray}
& & \dim _M(L^2(M)) = 1 , \label{eq-for01a} \\[1mm]
& & \dim _M(H)\cdot\dim _{M'}(H) = 1 , \label{eq-for01b} \\[1mm]
& & {\rm If} \ \{\pi , H\} \ {\rm and} \ \{\pi ', H'\} \ {\rm are \ faithful \ nondegenerate} \ W^*{\rm representations} \label{eq-for01c}\\[1mm]
& & {\rm of} \ M, \ {\rm then} \ \dim _M(H)=\dim _M(H') \ {\rm if \ and \ only \ if} \ \{\pi ,H\}\cong\{\pi ',H'\}, \nonumber \\[1mm]
& & {\rm i.e. \ if \ these \ W^*-representations \ are \ spatially \ *-isomorphic.} \nonumber \\[1mm]
& & {\rm If} \ \{\pi _i , H_i\}_{i\geq 1} \ {\rm is \ a \ sequence \ of \ faithful \ nondegenerate} \ W^*{\rm representa}\verb!-! \label{eq-for01d}\\[1mm]
& & {\rm tions \ of} \ M, \ {\rm then} \quad \dim _M\bigl(\sum _i H_i\bigr) = \sum _i\dim _M(H_i), \nonumber \\[1mm]
& & {\rm If} \ \{\pi , H\} \ {\rm is \ a \ faithful \ nondegenerate} \ W^*{\rm representation \ of} \ M, \ {\rm then}\label{eq-for01e}\\[1mm]
& & \pi(M)' \ {\rm is \ finite} \ {\rm if \ and \ only \ if} \ \dim _M(H)<\infty . \nonumber\\[1mm]
& & \dim _M(H)\geq 1 \ ({\rm resp.} \leq 1) \ {\rm if \ and \ only \ if} \  M \ {\rm admits \ a \ separating \
    (resp. \ cyclic)} \label{eq-for01f}\\[1mm]
& & {\rm vector}.\nonumber
\end{eqnarray}
\end{subequations}
}

\vspace{0.3cm}

We are now in a position to give the definition of the coupling constant for real
finite factors. Let us first prove an auxiliary Lemma.

\begin{lemma} \label{lemma001}
If $H$ is a real Hilbert space and $R\subset B(H)$ is a real
$W^*$-algebra, then $R' + iR' = (R + iR)'$, where
the latter commutant is taken in $B(H+iH)$.
\end{lemma}

\noindent {\bf Proof}. A straightforward calculation shows that
$R'+iR'\subset (R+iR)'$. \ Since $B(H+iH)= B(H)+iB(H)$ (see
\cite[Proposition 1.1.11]{Li2}), for each $a'\in (R+iR)'$ there
exist $x',y'\in B(H)$ such that $a'=x'+iy'$. \ Since $a'b=ba'$
for all $b=x+iy\in R+iR$, we have that $x'x-y'y = xx'-yy'$ and $x'x+y'y
= xx'+yy'$. Hence $x'x=xx'$ and $y'y=yy'$, i.e. $x',y'\in R'$.
Therefore $a'\in R'+iR'$. \ $\Box$\\

Now, let $M$ be a finite factor and let $\alpha$ be an involutive
*-antiautomor\-phism of $M$. If \ $\{\pi , H\}$ \ is a faithful
nondegenerate W*-repre\-sentation of $M$, and \ $\pi \circ \alpha =
{\tilde{\alpha}} \circ \pi$ \ for an involutive *-antiautomorphism
${\tilde{\alpha}}$ of $\pi(M)$, then by Lemma \ref{lemma001} we have \
$(\pi(M),{\tilde{\alpha}})' + i(\pi(M),{\tilde{\alpha}})' =
\pi(M)'$. Since the von Neumann algebra $\pi(M)'$ is semi-finite, the real factor
$(\pi(M),{\tilde{\alpha}})'$ is also semi-finite. Thus there exists
a unique (up to multiplication by a positive constant)
faithful normal semi-finite ${\tilde{\alpha}}$-invariant trace on
$\pi(M)'_+$. We define a {\it natural} ${\tilde{\alpha}}$-{\it invariant trace} on $\pi(M)'_+$ as follows.\\

{\bf (i)} If \ $\{\pi , H\} = \{\lambda\otimes \e ,  L^2(M)\otimes
K\}$, where $K$ is a countably infinite dimensional Hilbert space,
then the von Neumann algebra \ $(\lambda(M) \otimes \e _K)' = J\lambda(M)J
{\overline{\otimes}} B(K)$ \ is infinite and for the real factor \ $(M,\alpha)$ \ we
have
$$
\{\pi |_{(M,\alpha)}, H_r\} = \{\lambda _r\otimes \e ,
L^2(M,\alpha)\otimes K_r\} ,
$$
$$
(\lambda _r(M,\alpha) \otimes \e _{K_r})' = J\lambda _r(M,\alpha)J
{\overline{\otimes}} B(K_r)
$$
(further, for the sake of convenience, we shall write $\pi$, instead of $\pi |_{(M,\alpha)}$).
Pick an orthogonal normalized basis $\{e_i\}_{i\in
\Lambda}$ of $K$, where \ $|\Lambda |=\dim _{\mathbb{C}} K$. Then
each element \ $t'\in (\lambda(M) \otimes \e _K)'$ can be uniquely represented as
\ $t'=\bigl(J\lambda(x_{ij})J\bigr)$, \ where $x_{ij}\in M$, for all $i,j$.
If \ $t'\in (\lambda _r(M,\alpha) \otimes \e _{K_r})'$, i.e. \ ${\tilde{\alpha}}'(t') = (t')^*$,
\ then it is not difficult to show that $x_{ij}\in (M,\alpha)$ \
(i.e. \ $\alpha(x_{ij}) = x_{ij}^*$), for all $i,j$.
Define the natural trace as follows (see \cite[17.1.4 (i)]{Li1}):\\

{\tt Tr}$_{L^2(M)\otimes K}' (t') = \ \sum\limits _{i\in \Lambda}
\tau(x_{ii}), \quad \
t'=\bigl(J\lambda(x_{ij})J\bigr) \in (\lambda(M) \otimes \e _K)'_+$,\\

\noindent where $\tau$ is the unique faithful normal (and hence
$\alpha$-invariant) tracial state on $M$. It is easy to show that the
definition of \ {\tt Tr}$_{L^2(M)\otimes K}'$ \ is independent of
the choice of $\{e_i\}$ and \ {\tt Tr}$_{L^2(M)\otimes K}'$ \ is a
faithful semi-finite normal trace on \ $(\lambda(M) \otimes \e
_K)'_+$. Moreover, since\\[-1mm]

\noindent ({\tt Tr}$_{L^2(M)\otimes K}' \circ
{\tilde{\alpha}}')(t')=${\tt Tr}$_{L^2(M)\otimes
K}'({\tilde{\alpha}}'(t'))=$
{\tt Tr}$_{L^2(M)\otimes K}' ((t')^*)$\\

\qquad $= \sum\limits _i \tau(x_{ii}^*) = \ \sum\limits _i
\tau(\alpha(x_{ii})) = \ \sum\limits _i (\tau \circ \alpha)(x_{ii})
= \sum\limits _i \tau(x_{ii})$\\

\qquad = {\tt Tr}$_{L^2(M)\otimes K}' (t')$,\\

\noindent we have that \ {\tt Tr}$_{L^2(M)\otimes K}'$ \ is
${\tilde{\alpha}}'$-invariant. Therefore,
for {\tt Tr}$_{L^2(M,\alpha)\otimes K_r}'$ defined as follows\\

{\tt Tr}$_{L^2(M,\alpha)\otimes K_r}'(t')=\sum\limits _i
\tau(x_{ii})$, \
$t'=\bigl(J\lambda(x_{ij})J\bigr)\in (\lambda _r(M,\alpha)\otimes \e _{K_r})'_+$,\\

\noindent we have

\hspace{2.5cm} {\tt Tr}$_{L^2(M)\otimes K}' \Bigl |_{(\lambda
_r(M,\alpha)\otimes \e _{K_r})'} \Bigr. \ =$
\ {\tt Tr}$_{L^2(M,\alpha)\otimes K_r}'$ .\\[2mm]

{\bf (ii)} For a general faithful nondegenerate W$^*$-representation
$\{\pi , H\}$ of $M$ with $\pi \circ \alpha = {\tilde{\alpha}} \circ
\pi$ \ by Theorem \ref{teorem005} there are a projection $p'\in
(\lambda _r(M,\alpha)\otimes \e _{K_r})'$ \ and a unitary \ $u :
H_r\to p'(L^2(M,\alpha)\otimes K_r)$ \ such that
$$
u\pi(x)u^* = (\lambda(x)\otimes \e _K)p' , \qquad \ x\in M ,
$$
where $K_r$ is a real Hilbert space and $K=K_r+iK_r$.
Then we define the natural trace as follows (see \cite[17.1.4 (ii)]{Li1}):\\

\hspace{2.2cm} {\tt Tr}$_H' (t') =$ \ {\tt Tr}$_{L^2(M)\otimes K}' (ut'u^*), \qquad \ t'\in \pi(M)'_+$ .\\

\noindent The definition of {\tt Tr}$_H'$ is independent on the
choice of $u$ and $p'$, and {\tt Tr}$_H'$  is a faithful normal
trace on $\pi(M)'_+$. Since ${\tilde{\alpha}}'(u)=u^*$, we have that {\tt
Tr}$_H'$  is ${\tilde{\alpha}}'$-invariant. Therefore,
for {\tt Tr}$_{L^2(M,\alpha)\otimes K_r}'$ defined as\\

\hspace{0.75cm}  {\tt Tr}$_{L^2(M,\alpha)\otimes K_r}' (ut'u^*) =$
{\tt Tr}$_{L^2(M)\otimes K}' (ut'u^*) ,
\qquad \ t'\in (\pi(M),{\tilde{\alpha}})'_+$\\

\noindent we have

\hspace{2.75cm} {\tt Tr}$_H' \Bigl|_{(\pi(M),{\tilde{\alpha}})'} \Bigr. \ =$ \ {\tt Tr}$_{L^2(M,\alpha)\otimes K_r}'$ .\\[2mm]

\noindent If {\tt Tr}$_{H_r}'$ denotes {\tt Tr}$_{L^2(M,\alpha)\otimes K_r}'$, then we have

\hspace{2.75cm} \ {\tt Tr}$_{H_r}'$ = {\tt Tr}$_H' \Bigl|_{(\pi(M),{\tilde{\alpha}})'} \Bigr.$ .\\[2mm]

\noindent Thus, \ {\tt Tr}$_{H}' \Bigl|_{(\pi(M),{\tilde{\alpha}})'}
\Bigr.$ is a faithful normal semi-finite trace on
$(\pi(M),{\tilde{\alpha}})'$.


\begin{definition} \label{definition01}
Let $M$ be a finite factor and let $\alpha$ be an involutive
*-antiautomorphism of $M$. Suppose that \ $\{\pi , H\}$ \ is a faithful
nondegenerate W*-representation of $M$, and \ $\pi \circ \alpha =
{\tilde{\alpha}} \circ \pi$ \ for an involutive *-antiautomorphism \
${\tilde{\alpha}}$ of $\pi(M)$. Then\\[-1mm]

\hspace{3.5cm} $\dim _{(M,\alpha)}(H_r) =$ {\tt Tr}$_{H_r}'(\e)$\\

\noindent is called the {\bf coupling constant} between
\ $(\pi(M),{\tilde{\alpha}})$ and $(\pi(M),{\tilde{\alpha}})'$ \ relative to $H_r$.
\end{definition}

\vspace{0.25cm}

{\bf (iii)} Now, in the case where \ $\{\pi , H\} = \{\lambda\otimes \e ,
L^2(M)\otimes K\}$ we choose another basis. Namely, pick a real
orthogonal normalized basis $\{f_i\}_{i\in \Lambda '}$ of $K$, where
\ $|\Lambda '|=\dim _{\mathbb{R}} K$. Then
each element \ $t'\in (\lambda(M) \otimes \e _K)'$
can be uniquely represented as \
$t'=\bigl(J\lambda(x_{ij})J\bigr)$, \ where $x_{ij}\in M$, for all
$i,j$. If \ $t'\in (\lambda _r(M,\alpha) \otimes \e _{K_r})'$, i.e. \
${\tilde{\alpha}}'(t') = (t')^*$, \ then it is not difficult to show that
\ $x_{ij}\in (M,\alpha)$ \ (i.e. \ $\alpha(x_{ij}) = x_{ij}^*$),
for all $i,j$. We set\\

{\tt tr}$_{L^2(M)\otimes K}' (t') = \ \sum\limits _{i\in \Lambda '}
\tau(x_{ii}), \quad \ t'=\bigl(J\lambda(x_{ij})J\bigr) \in (\lambda(M) \otimes \e _K)'_+$ .\\

\noindent Clearly, \ {\tt tr}$_{L^2(M)\otimes K}'$ \ is also a
faithful normal semi-finite trace on \ $(\lambda(M) \otimes
\e _K)'_+$. Moreover {\tt tr}$_{L^2(M)\otimes K}'$ \ is
${\tilde{\alpha}}'$-invariant. Similarly, we can show that the
definition of \ {\tt tr}$_{L^2(M)\otimes K}'$ \ does not depend on
the choice of $\{f_i\}$.\\[1mm]
In the case where \ $u\pi(x)u^* = (\lambda(x)\otimes \e _K)p'$ \ ($x\in M$) \ we put\\[-2mm]

\hspace{2.2cm} {\tt tr}$_H' (t') =$ \ {\tt tr}$_{L^2(M)\otimes K}' (ut'u^*), \qquad \ t'\in \pi(M)'_+$ .\\

\noindent The definition of {\tt tr}$_H'$ is also independent on the
choice of $u$ and $p'$, and the trace {\tt tr}$_H'$  is ${\tilde{\alpha}}'$-invariant.


Let \ $\{\pi , H\}$ \ be a faithful nondegenerate W*-representation of $M$ \ with \ $\pi \circ \alpha =
{\tilde{\alpha}} \circ \pi$ \ for an involutive *-antiautomorphism \ ${\tilde{\alpha}}$ of $\pi(M)$.

\vspace{0.5cm}

\begin{definition} \label{definition02}
The number \ $\dim _{(M,\alpha)}(H) =$ {\tt tr}$_{H}'(\e)$ \ is called the {\bf coupling constant} between
$(\pi(M),{\tilde{\alpha}})$ and $(\pi(M),{\tilde{\alpha}})'$ \ relative to $H$.
\end{definition}

One has the following relations between \ $\dim _{(M,\alpha)}(H_r)$, \ $\dim _{(M,\alpha)}(H)$ \ and \
$\dim _M(H)$.

\begin{theorem} \label{teorem0061}
$$
\dim _M(H) = \dim _{(M,\alpha)}(H_r) = \frac{1}{2} \dim
_{(M,\alpha)}(H)
$$
\end{theorem}

\noindent The proof of this theorem is obvious.  \ $\Box$ \\[-2mm]

Let's consider some properties of the coupling constants.

\vspace{0.35cm}

\begin{proposition}  \label{proposition3}
Let $M$ $(\subset B(H))$ be a finite factor and let $\alpha$ be an involutive *-antiautomorphism of $M$. Then\\[1mm]
{\rm (i)} \ $\dim _{(M,\alpha)}(L^2(M))=2$ \ and \ $\dim _{(M,\alpha)}(L^2(M,\alpha))=1$.\\[2mm]
{\rm (ii)} \ $\dim _{(M,\alpha)}(H)\cdot\dim _{(M,\alpha)'}(H)=4$ and $\dim _{(M,\alpha)}(H_r)\cdot\dim _{(M,\alpha)'}(H_r)=1$\\[2mm]
{\rm (iii)} If $\{\pi , H\}$ and $\{\pi ', H'\}$ are $\alpha$-invariant faithful nondegenerate W*-representations of $M$, then
\ $\dim _{(M,\alpha)}(H)=\dim _{(M,\alpha)}(H ')$ \ if and only if \ $\{\pi , H\}$ \ and \ $\{\pi ', H'\}$ \
are spatially *-isomorphic via a unitary operator $w$ with \ $\pi(\alpha(w)) = {\hat{\alpha}} \pi '(w) = \pi '(w)^*$ ;\\[2mm]
{\rm (iv)} If $\{\pi _i , H_i\}_{i\geq 1}$ is a sequence of
$\alpha$-invariant faithful nondegenerate W*-representations of $M$,
then
\ $\dim _{(M,\alpha)}(\sum\limits _iH_i)=\sum\limits _i\dim _{(M,\alpha)}(H_i)$;\\[2mm]
{\rm (v)} If $\{\pi , H\}$ is an $\alpha$-invariant faithful nondegenerate W*-representation of $M$, then
the following conditions are equivalent:

{\rm a)} real von Neumann algebra \ $(\pi(M),{\tilde{\alpha}})'$ is finite;

{\rm b)} the trace $Tr'_{H_r}$ is finite;

{\rm c)} $\dim _{(M,\alpha)}(H) < \infty$ .\\[2mm]
{\rm (vi)} $\dim _{(M,\alpha)}(H)\geq 2$ \ (resp. $\leq 2$) \
if and only if \ $(M,\alpha)$ admits a separating (resp. \
cyclic) vector.
\end{proposition}

\vspace{0.5cm}

\noindent {\bf Proof}. The property \eqref{eq-for01a} and Theorem \ref{teorem0061} imply the proof
of (i). From \eqref{eq-for01b} and Theorem \ref{teorem0061} we obtain
the proof of (ii). The equivalence of the conditions \ $\dim _{(M,\alpha)}(H)=\dim _{(M,\alpha)}(H ')$ \
and \ $\dim _M(H)=\dim _M(H ')$ \ follows from Theorem \ref{teorem0061}.
The equivalence of the conditions \ $\dim _M(H)=\dim _M(H ')$ \ and \ $\{\pi , H\}  \cong \{\pi ', H'\}$ \
follows from \eqref{eq-for01c}.
By Theorem \ref{teorem005} in this case there exists a unitary operator $w$ with
$\pi(\alpha(w))={\hat{\alpha}} \pi '(w) = \pi '(w)^*$ \ which implements this spatial *-isomorphism
\ $\{\pi , H\}\cong \{\pi ', H'\}$, \ what is required to be proved for (iii).
From \eqref{eq-for01d} we obtain the proof of (iv).
From the theory of real W*-algebras we know that \
the real von Neumann algebra $(\pi(M),{\tilde{\alpha}})'$ is finite if and only if \
the von Neumann algebra $\pi(M)'$ is finite (see \cite{ARU}).
Then by \eqref{eq-for01e} we obtain the proof of (v).
It is easy to see that $(M,\alpha)$ admits a separating (respectively, cyclic)
vector if and only if $M$ admits a separating (respectively, cyclic)
vector. Then from \eqref{eq-for01f} and Theorem \ref{teorem0061} we
obtain the proof of (vi). \ $\Box$

\begin{proposition}  \label{proposition4}
If $\Re$ is a finite real factor on a real Hilbert space $H$ with the
finite commutant $\Re '$, and $\tau$, $\tau '$ are the unique
faithful normal tracial states on $\Re$ and $\Re '$ respectively,
then for any $\xi (\not=0) \in H$ the number \ $\displaystyle c_\Re
= \frac{\tau(e_\xi)}{\tau '(e'_\xi)}$ \ is independent of the choice
of $\xi$. Moreover, we have \ $c_\Re = \dim _\Re(H)$, where $e_\xi$
and $e'_\xi$ are the cyclic projections from $H$ onto
${\overline{\Re '\xi}}$ \ and \ ${\overline{\Re \xi}}$ \ respectively.
\end{proposition}

\noindent {\bf Proof}. We extend $\tau$ and $\tau '$ to $\Re +
i\Re$ and $\Re '+ i\Re '$, respectively, by the linearity as
${\overline{\tau}}(a+ib)=\tau(a)+i\tau(b)$ and
${\overline{\tau}}'(a'+ib')=\tau '(a')+i\tau '(b')$. It is
obvious, that for the cyclic projections ${\overline{e}}_\xi$ and
${\overline{e}}'_\xi$ from $H_c=H+iH$ onto ${\overline{\Re '\xi}} +
i{\overline{\Re '\xi}}$, \ ${\overline{\Re \xi}}+i{\overline{\Re
\xi}}$, \ respectively, we have $\displaystyle
\frac{{\overline{\tau}}({\overline{e}}_\xi))}{{\overline{\tau}}'({\overline{e}}'_\xi)}
= \frac{\tau(e_\xi)}{\tau '(e'_\xi)}$. \ Since $\displaystyle
\frac{{\overline{\tau}}({\overline{e}}_\xi))}{{\overline{\tau}}'({\overline{e}}'_\xi)}$
\ is independent on the choice of $\xi$, and this number is equal
to $\dim _{\Re +i\Re}(H+iH)$, \ we obtain that $\dim _\Re(H)= c_\Re$. \
$\Box$

\vspace{0.75cm}

\begin{example} \label{example1}
{\rm Let $M$ be a factor of type I$_n$ on an $m$-dimensional Hilbert
space $H$ ($n,m<\infty$), and let $\alpha$ be an involutive
*-antiautomorphism of $M$. It is known that $\dim _M(H)=m/n^2$, i.e
$\dim _M(H)=\dim(H)/\dim(M)$. Then by Theorem \ref{teorem0061} we
have $\dim _{(M,\alpha)}(H)=2m/n^2$, hence, since
$\dim(M,\alpha)=n^2$, \ one has \ $\dim _{(M,\alpha)}(H) = \frac{2
\dim(H)}{\dim(M)} = \dim _{\mathbb{R}}(H)/\dim(M,\alpha)$. }
\end{example}
From the theory of real W*-algebras we know (see \cite{ARU,G,R}),
that if $n$ is an odd number, then $M$ possesses exactly one conjugacy
class of involutive *-antiautomorphism, and if $n$ is an even number,
then there are two conjugacy class of involutive *-antiautomorphism
of $M$. Namely, in the first case $(M,\alpha) \cong
M_n(\mathbb{R})$, and in the second case we have also the possibility
$(M,\alpha) \cong M_{n/2}(\mathbb{H})$, where $\mathbb{H}$ is the quaternion algebra. Thus\\[-1mm]

$\displaystyle \dim _{M_n(\mathbb{R})}(H) = \frac{\dim
_{\mathbb{R}}(H)}{n^2} =
\frac{\dim _{\mathbb{R}}(H)}{\dim(M_n(\mathbb{R}))}$ ,\\[2mm]

$\displaystyle \dim _{M_{n_{/2}}(\mathbb{H})}(H) = \frac{4\dim
_{\mathbb{H}}(H)}{n^2} = \frac{\dim
_{\mathbb{H}}(H)}{\bigl(\frac{n}{2}\bigr)^2} = \frac{\dim
_{\mathbb{H}}(H)}{\dim _{\mathbb{H}}(M_{n_{/2}}(\mathbb{H}))}$ .


\section{The index of subfactors in finite real factors}

\begin{definition} \label{definition03}
Let $M\subset B(H)$ be a finite factor. Consider a subfactor $N\subset M$ and let
$\alpha$ be an involutive *-antiautomorphism of $M$ with $\alpha(N)\subset N$.
Consider the real factors $\Re=(M, \alpha)$ and $Q=(N, \alpha)$. The {\bf index}
of $Q$ in $\Re$, denoted by $[\Re : Q]$ or by $[(M, \alpha) : (N, \alpha)]$,
is defined as \ $\dim _Q(L^2(\Re))$.
\end{definition}

Between real and complex indices there is the following relation

\begin{theorem} \label{teorem041}
We have \ $[(M,\alpha):(N,\alpha)]=[M,N]$, i.e. \ $[\Re : Q]=[\Re +
i\Re : Q+iQ]$.
\end{theorem}

\noindent {\bf Proof.} By Theorem \ref{teorem0061} we have \
$[(M,\alpha):(N,\alpha)] = \dim _{(N,\alpha)}(L^2(M,\alpha)) =
\frac{1}{2} \dim _{(N,\alpha)}(L^2(M)) =
\frac{1}{2} \cdot 2 \dim _N(L^2(M)) = [M,N]$ \ $\Box$\\

We recall some properties of the complex index (\cite{J}, \cite{K1},
\cite[Ch. 17]{Li1}): if $M$ is a finite factor and $N$ is a subfactor of $M$
then

\vspace{-0.3cm}

{\small
\begin{subequations} \label{eq-for07}
\begin{eqnarray}
& & [M:N] = \ \dim _N(H)/\dim _M(H), \label{eq-for07a} \\[1mm]
& & [M:N] \geq \ 1, \ {\rm in \ particular}, \ [M:M]=1, \label{eq-for07b} \\[1mm]
& & [M:N]=[N':M'], \label{eq-for07c}\\[1mm]
& & {\rm If} \ P \ {\rm is \ a \ subfactor \ of} \ N, \ {\rm then} \ [M:P]=[M:N]\cdot[N:P], \label{eq-for07d}\\[1mm]
& & {\rm If} \ P \ {\rm is \ a \ subfactor \ of} \ N, \ [M : P]<\infty \ {\rm and} \ [M:P]=[M:N], \label{eq-for07e}\\[1mm]
& & {\rm then} \ N=P, \nonumber \\[1mm]
& & {\rm If} \ M_i \ {\rm is \ a \ finite \ factor \ and} \ N_i \ {\rm is \ a \ subfactor \ of} \ M_i (i=1,2) \label{eq-for07f}\\[1mm]
& & {\rm then} \
[M_1{\overline\otimes}M_2:N_1{\overline\otimes}N_2]=[M_1:N_1]\cdot
[M_2:N_2] \nonumber
\end{eqnarray}
\end{subequations}
}

\vspace{0.3cm}

Similarly to the complex case, using Theorem \ref{teorem041} we can
prove the following pro\-perties of the real index

\begin{theorem} \label{teorem042}
Let $H_r$ be a real Hilbert space. Suppose that \ $\Re\subset B(H_r)$ is a finite
real factor and $Q\subset \Re$ is a real subfactor. Let $M=\Re +i\Re$
be the enveloping complex factor for $\Re$ and let $\alpha$ be the involutive
*-antiautomorphism of $M$ which generates $\Re$, i.e.
$\Re=(M,\alpha)$ (in this case it is clear that $Q=(N,\alpha)$, where $N=Q+iQ$).
Then\\[3mm]
{\rm (i)} \ $[(M,\alpha):(N,\alpha)] \ = \ \displaystyle \frac{\dim
_{(N,\alpha)}(H_r)}{\dim _{(M,\alpha)}(H_r)}$ , i.e.
\ $[\Re : Q] \ = \ \displaystyle \frac{\dim _Q(H_r)}{\dim _\Re(H_r)}$ .\\[3mm]
{\rm (ii)} \ $[(M,\alpha):(N,\alpha)] \ \geq \ 1$ , i.e. \ $[\Re :
Q] \ \geq \ 1$ . In particular,
\ $[(M,\alpha):(M,\alpha)] = [\Re : \Re] \ = \ 1$ .\\[3mm]
{\rm (iii)} \ $[(M,\alpha):(N,\alpha)] = [(N,\alpha)' : (M,\alpha)']$, \ i.e. \ $[\Re : Q] = [Q' : \Re ']$.\\[3mm]
{\rm (iv)} \ If $Q_1$ is a real subfactor of $Q$, then \
$[(M,\alpha):(N_1,\alpha)] = [(M,\alpha) : (N,\alpha)] \cdot
[(N,\alpha) : (N_1,\alpha)]$, \
i.e. \ $[\Re : Q_1] = [\Re : Q] \cdot [Q : Q_1]$, \ where $N_1=Q_1+iQ_1$.\\[3mm]
{\rm (v)} \ If $Q_1$ is a real subfactor of $Q$, \ $[\Re :
Q_1]<\infty$ and \ $[(M,\alpha):(N_1,\alpha)] = [(M,\alpha) :
(N,\alpha)]$, \ then $N=N_1$, \ i.e.
\ if \ $[\Re : Q_1] = [\Re : Q]$, \ then $Q=Q_1$, \ where $N_1=Q_1+iQ_1$.\\[4mm]
{\rm (vi)} Let \ $\Re _i$ \ be a finite real factor, and let $Q_i$ be a real
subfactor of $\Re _i$, $i=1,2$. If $M_i=\Re _i +i\Re _i$ and
$N_i=Q_1+iQ_i$ are the enveloping complex factors for $\Re _i$ and $Q_i$, respectively,
then let $\alpha _i$ denote the involutive *-antiautomorphism of $M_i$,
which gene\-rates $\Re _i$, i.e. $\Re _i=(M_i,\alpha _i)$, $i=1,2$. Then\\[2mm]
$[(M_1,\alpha _1){\overline\otimes}(M_2,\alpha _2) : (N_1,\alpha
_1){\overline\otimes}(N_2,\alpha _2)] = [(M_1,\alpha _1):(N_1,\alpha
_1)] \cdot [(M_2,\alpha _2):(N_2,\alpha _2)]$, \ i.e.
\ $[\Re _1{\overline\otimes} \Re _2 : Q_1{\overline\otimes}Q_2] = [(\Re _1 : Q_1] \cdot [\Re _2 : Q_2]$ .\\[2mm]
\end{theorem}

\noindent {\bf Proof.} (i). By Theorem \ref{teorem041} and the property
(\ref{eq-for07a}) we have $[(M,\alpha):(N,\alpha)] =$\\

$=[M:N] =
\displaystyle \frac{\dim _N(H_r+iH_r)}{\dim _M(H_r+iH_r)} =
\displaystyle \frac{\frac{1}{2}\dim
_{(N,\alpha)}(H_r+iH_r)}{\frac{1}{2}\dim _{(M,\alpha)}(H_r+iH_r)} =
\displaystyle \frac{\dim _{(N,\alpha)}(H_r)}{\dim _{(M,\alpha)}(H_r)}$.\\[1mm]

\noindent (ii). \ By the Theorem \ref{teorem041} and
(\ref{eq-for07b}) we have
$[(M,\alpha):(N,\alpha)] = [M:N] \geq 1$ and $[(M,\alpha):(M,\alpha)] = [M:M] = 1$.\\[-1mm]

\noindent (iii). \ As above, let $\alpha '$ be the involutive
*-antiautomorphism of $M'$, which generates $\Re '$, i.e. $\Re
'=(M',\alpha ')$. Then $(M,\alpha)' = (M',\alpha ')$. Similarly we
have $(N,\alpha)' = (N',\alpha ')$. Hence by Theorem
\ref{teorem041} and the property (\ref{eq-for07c}) we have
$[(M,\alpha):(N,\alpha)] = [M:N] = [N':M'] = [(N',\alpha ') : (M',\alpha ')] = [(N,\alpha)' : (M,\alpha)']$.\\[-1mm]

\noindent (iv).\ By Theorem \ref{teorem041} and the property
(\ref{eq-for07d}) we have \ $[(M,\alpha):(N_1,\alpha)] = [M:N_1] =
[M:N]\cdot[N:N_1] =
[(M,\alpha) : (N,\alpha)] \cdot [(N,\alpha) : (N_1,\alpha)]$.\\[-1mm]

\noindent (v).\ If \ $[(M,\alpha):(N_1,\alpha)] = [(M,\alpha) :
(N,\alpha)]$, by Theorem \ref{teorem041} we have \
$[M:N_1]=[(M,\alpha):(N_1,\alpha)] = [(M,\alpha) : (N,\alpha)] =
[M:N]$, i.e.
$[M:N_1]=[M:N]$. Then by (\ref{eq-for07e}) we obtain $N=N_1$, i.e. $Q=Q_1$.\\[-1mm]

\noindent (vi).\ According to \cite[page 69]{Li2} we have $\Re
_1{\overline\otimes} \Re _2 + i \Re _1{\overline\otimes} \Re _2 =
(\Re _1 + i\Re _1 ){\overline\otimes}(\Re _2 + i\Re _2 )$, i.e.
$(M_1,\alpha _1){\overline\otimes} (M_2,\alpha _2) + i (M_1,\alpha
_1){\overline\otimes} (M_2,\alpha _2)  = M_1{\overline\otimes}M_2$.
Similarly, we have \ $(N_1,\alpha _1){\overline\otimes} (N_2,\alpha
_2) + i (N_1,\alpha _1){\overline\otimes} (N_2,\alpha _2)  =
N_1{\overline\otimes}N_2$.
Then by Theorem \ref{teorem041} and the pro\-perty (\ref{eq-for07f}) we get\\[2mm]
$[(M_1,\alpha _1){\overline\otimes}(M_2,\alpha _2) : (N_1,\alpha _1){\overline\otimes}(N_2,\alpha _2)] =$\\[-2mm]

\noindent $[(M_1,\alpha _1){\overline\otimes}(M_2,\alpha _2) + i
(M_1,\alpha _1){\overline\otimes}(M_2,\alpha _2) \ : \
(N_1,\alpha _1){\overline\otimes}(N_2,\alpha _2) +$\\[-2mm]

\noindent $i (N_1,\alpha _1){\overline\otimes}(N_2,\alpha _2)] =
[((M_1,\alpha _1) + i (M_1,\alpha _1)) {\overline\otimes} ((M_2,\alpha _2) + i (M_2,\alpha _2)):$\\[-2mm]

\noindent $((N_1,\alpha _1) + i (N_1,\alpha _1)) {\overline\otimes}
((N_2,\alpha _2) + i (N_2,\alpha _2))]=
[M_1{\overline\otimes}M_2 : N_1{\overline\otimes}N_2]=$\\[-2mm]

\noindent  $[M_1:N_1]\cdot [M_2:N_2] = [(M_1,\alpha _1):(N_1,\alpha _1)] \cdot [(M_2,\alpha _2):(N_2,\alpha _2)]$.\\[2mm]
Thus, we have
$[\Re _1{\overline\otimes} \Re _2 : Q_1{\overline\otimes}Q_2] = [(\Re _1 : Q_1] \cdot [\Re _2 : Q_2]$ . \ $\Box$\\

As it was noted in the introduction, V.Jones in \cite{J} has proved a theorem on the values of index for subfactors.
Let us recall this theorem

\vspace{0.5cm}

\begin{theorem}[\cite{J}, Theorem 4.3.1] \label{teorem043}
Let $M$ be a finite factor, and let $N$ be a subfactor of $M$ with
$[M:N]<\infty$. Then one has either \ $\displaystyle [M:N]=4\cos
^2\frac{\pi}{q}$ \ for some integer \ $q\geq 3$ \ or \ $[M:N]\geq
4$.
\end{theorem}

From Theorems \ref{teorem041} and \ref{teorem043} we obtain the following real version
of the above theorem.

\begin{theorem} \label{teorem044}
Let $M$ be a finite factor and let $N$ be a subfactor of $M$ with
$[M:N]<\infty$. Given be an involutive *-antiautomorphism $\alpha$ of
$M$ with $\alpha(N)\subset N$, put $\Re=(M,\alpha)$, $Q=(N,\alpha)$.
Then one has either \ $\displaystyle
[(M,\alpha):(N,\alpha)]=4\cos ^2\frac{\pi}{q}$ \ for some integer \
$q\geq 3$ \ or \ $[(M,\alpha):(N,\alpha)]\geq 4$, i.e. \
$\displaystyle [\Re :Q]=4\cos ^2\frac{\pi}{q}$ \ for some integer \
$q\geq 3$ \ or \ $[\Re :Q]\geq 4$.
\end{theorem}

\vspace{2.5cm}

\textbf{Ackowledgments.} \emph{The second, the third and the fourth named authors would like to acknowledge the hospitality
of the "Institut f$\ddot{u}$r Angewandte Mathematik",
Universit$\ddot{a}$t Bonn (Germany). This work is supported in
part by the DFG 436 USB 113/10/0-1 project (Germany) and the
Fundamental Research Foundation of the Uzbekistan Academy of
Sciences.}

\end{document}